\documentclass[12pt]{article}
\advance\voffset by -0.7cm
\advance\hoffset by -1.7cm
\input amssym.def

\def\N{{\Bbb N}}     
\def\Z{{\Bbb Z}}
\def\Qq{{\Bbb Q}}    
\def\R{{\Bbb R}}
\def\C{{\Bbb C}}     
\def\H{{\Bbb H}}    
\def\D{{\Bbb D}}
\def\SP{\mbox{{\rm{Sp}}}(2n,\Z)}
\def\SPQ{\mbox{{\rm{Sp}}}(2n,\Qq)}
\def\SPR{\mbox{{\rm{Sp}}}(2n,\R)}
\def\GL#1{\mbox{GL}(#1,\C)}
\def\P{{\cal P}}
\def\Q{{\cal Q}}
\def\Hq{{\cal H}}
\def\bin#1#2{{\left( #1\atop #2 \right)}}
\def\Sym{\mbox{\rm Sym}}
\def\min{\mbox{\rm min}}
\def\tr{{\mbox{\rm tr}}}
\def\tfrac#1#2{{\scriptstyle{\frac{#1}{#2}}}}
\def\sq{\hbox{\rlap{$\sqcap$}$\sqcup$}}
\def\qed{\ifmmode\sq\else{\unskip\nobreak\hfil
         \penalty50\hskip1em\null\nobreak\hfil\sq
         \parfillskip=0pt\finalhyphendemerits=0\endgraf}\fi}
\def\id{{\mathchoice {\rm 1\mskip-4mu l} {\rm 1\mskip-4mu l}
        {\rm 1\mskip-4.5mu l} {\rm 1\mskip-5mu l}}}

\newtheorem{thm}{Theorem}[section]
\newtheorem{lem}[thm]{Lemma}
\newtheorem{prop}[thm]{Proposition}
\newtheorem{df}[thm]{Definition}

\begin{document}
\title{ Rankin-Cohen Type Differential Operators \\
        for Siegel Modular Forms}
\author{ Wolfgang Eholzer
         \thanks{Supported by the EPSRC and partially by 
                  PPARC and  EPSRC (grant GR/J73322).}\\
 {\small Department of Applied Mathematics and Theoretical Physics,}\\
 {\small University of Cambridge, Silver Street,}\\{\small Cambridge
       CB3 9EW, U.K.}\\
 {\small e-mail: W.Eholzer@damtp.cam.ac.uk}
\and  
\\
        Tomoyoshi Ibukiyama  \\
{\small Department of Mathematics, Graduate School of Mathematics,}\\
{\small Osaka University, Machikaneyama 1-16,}\\
{\small Toyonaka, Osaka, 560 Japan}\\
{\small e-mail: ibukiyam@math.wani.osaka-u.ac.jp}
}
\date{}
\maketitle
\abstract{
Let $\H_n$ be the Siegel upper half space and let $F$ and $G$ be 
automorphic forms on $\H_n$ of weights $k$ and  $l$, respectively.
We give explicit examples of differential operators $D$ acting
on functions on $\H_n \times\H_n$ such that the restriction of  
$$ D(F(Z_{1})G(Z_{2})) $$
to $Z = Z_{1} = Z_{2} $ 
is again an automorphic form of weight $k + l + v$ on $\H_n$.
Since the elliptic case, {\it i.e.} $n=1$, has already been studied 
some time ago by R.\ Rankin and H.\ Cohen we call such differential
operators Rankin-Cohen type operators.

We also discuss a generalisation of Rankin-Cohen type  
operators to vector valued differential operators.
}

\vspace{1.0cm}
\begin{flushright}
DAMTP-97-26 \\
alg-geom/9703033
\end{flushright}

\section{Introduction}

In this paper  we are concerned with the explicit construction
of bilinear differential operators for Siegel modular forms mapping 
$M(\Gamma)_k\times M(\Gamma)_l$ to  $M(\Gamma)_{k+l+v}$
for all even non-negative integers $v$ and $\Gamma$
some discrete subgroup of $\SPR$ with finite co-volume. 
It has been shown in ref.\ \cite{Ibukiyama} that if the weights $k$ 
and $l$ are sufficiently large then there is a one-to-one 
correspondence between such bilinear differential operators and certain
invariant pluri-harmonic polynomials (for more details see \S\ref{basics}).
Therefore, one possibility for constructing such differential
operators is to construct the corresponding invariant pluri-harmonic 
polynomials.
As we will choose exactly this way for our construction one may view 
our paper also as an attempt to describe certain spaces of invariant 
pluri-harmonic polynomials explicitly.

Some results concerning the explicit construction of 
bilinear differential operators for Siegel modular 
forms --which we will call Rankin-Cohen type operators following
ref.\ \cite{Zagier}-- are already known in the literature: 
the genus one case has been considered 
R.\ Rankin \cite{Rankin} and H.\ Cohen \cite{Cohen} and, more recently, 
the genus two case by  Y.\ Choie and the first author~\cite{ChoieEholzer}.  
Both of these results are special cases of the construction 
presented in this paper. 
We are, however, not able to give closed explicit formulas 
for all Rankin-Cohen type operators, {\it i.e.} for
general values of $n$ and $v$. Instead we derive a system of 
recursion equations which is indeed very simple to solve 
for any {\it numerical} values of $n$ and $v$. In several special
cases, however, we obtain explicit closed formulas.
In addition we show that the image of the Rankin-Cohen type 
operators is, for $v>0$, contained in the space of cusp forms.
Finally, we discuss vector valued generalisations of the Rankin-Cohen
type operators, {\it i.e.} bilinear differential operators which map two
Siegel modular forms to a vector valued modular Siegel modular form
that transforms under a certain representation $\rho$ of $\GL{n}$.

This paper is organised as follows.
In section \ref{basics}, after reviewing some basic definitions 
and standard notations,  we recall the one-to-one correspondence 
between invariant pluri-harmonic polynomials and 
covariant differential operators. 
In \S\ref{main} we derive our main result in form  of a set of 
recursion relations which allow to determine certain invariant 
pluri-harmonic polynomials $Q_{n,v}$ (Theorem \ref{vgeneral}).
We also show that, for $v>0$, the bilinear differential operators 
associated to these polynomials map two Siegel modular forms to a 
Siegel cusp form. 
Section \ref{proofs} contains the proof of a uniqueness result 
in \S\ref{basics}  
as well as the proofs of three lemmas and our main theorem in 
\S\ref{main} and the result about cusp forms. 
We then discuss in \S\ref{examples} several special cases of our 
main result explicitly: the cases $v=2,4$ for general $n$ and the 
cases $n=1,2$ for general $v$.
Finally, we generalise Rankin-Cohen type differential operators to 
vector valued differential operators in section \ref{vector}. 
Here the case $n=2$ is treated explicitly and in some detail. 
We conclude in section \ref{conclusion} with some remarks and open questions.

\section{Siegel modular forms, pluri-harmonic polynomials 
         and differential operators}
\label{basics}

In this section we recall some standard notations 
and review a result of the second author on the relation
between invariant pluri-harmonic polynomials and 
differential operators (Theorem \ref{Ibuthm}).

Let $\H_n$ be the space of complex symmetric $n\times n$ matrices 
with positive definite imaginary part and define an action 
of $\SPR$ on functions $f:\H_n\to\C$ by
$$ (f|_M^k)(Z) = f(MZ) \det(CZ+D)^{-k} \qquad (M\in\SPR)$$
where $Z\in\H_n$, 
$M=\left( \begin{array}{cc} A & B \\ C & D \end{array}\right)$
with $n\times n$ matrices $A,B,C,D$, and, where
$MZ = (AZ+B)(CZ+D)^{-1}$.

We define Siegel modular forms on some discrete subgroup 
$\Gamma\subset \SPR$ of finite co-volume (finite co-volume 
means that the volume of $\SPR/\Gamma$ is finite).
\begin{df}
A holomorphic function $f:\H_n\to\C$ is called a Siegel modular 
form of non-negative weight $k$ on $\Gamma$ if 
$$ (f|_M^k)(Z) = f(Z)\qquad {\rm{for}}\ {\rm{all}}\  M\in\Gamma $$
(and $f(Z)$ is bounded at the cusps for $n=1$).
\end{df}
We denote the space of all Siegel modular forms of weight $k$ on
$\Gamma$ by $M(\Gamma)_k$. 

We review some notations concerning pluri-harmonic polynomials.
\begin{df} 
Let $P$ be a polynomial in the matrix variable  $X = (x_{r,s})\in M_{n,d}$
and define 
$$ \Delta_{i,j}(X) = \sum_{\nu=1}^{d} 
                   \frac{\partial^2}
                        {\partial  x_{i,\nu} \partial x_{j,\nu}}
   \qquad (1\le i,j\le n).
$$
Then the polynomial $P$ is called harmonic if 
$\sum_{i=1}^n \Delta_{i,i}(X) P = 0$ and $P$ is called 
pluri-harmonic if $\Delta_{i,j}(X) P = 0$ for all $1\le i,j\le n$. 
\end{df}
The group $\GL{n}\times O(d)$ acts on such polynomials  by
$P(X) \to P(A^t X B)$ ($A\in\GL{n}$, $B\in O(d)$).
We will be interested in pluri-harmonic polynomials
which are invariant under a subgroup of $O(d)$ and transform 
under a certain representation of $\GL{n}$.

A polynomial $P$  of a matrix variable $X \in M_{n,d}$ 
is called homogeneous of weight $v$ 
if $P(A^t X) = \det(A)^v\, P(X)$ for all $A\in\GL{n}$.
It was already pointed out in ref.\ \cite{KashiwaraVergne} 
that a polynomial $P(X)$ is pluri-harmonic if and only if 
$P(A X)$ is harmonic for all $A\in\GL{n}$.
Therefore, a homogeneous polynomial of some weight $v$ is 
pluri-harmonic if and only if it is harmonic.  

Let $d_i$ ($1\le i\le r$) be natural numbers such that $d_i\ge n$
and $d_1+\dots+d_r = d$. Define an embedding 
of $K=O(d_1)\times\dots\times O(d_r)$ into $O(d)$ by
$$ (B_1,\dots, B_r) \to 
  \left(\begin{array}{ccccc}
         B_1    &   0    & \dots & \dots & 0 \cr
         0      & B_2    &   0   & \dots & 0 \cr
         \vdots &   0    & \ddots&       & \vdots\cr
         \vdots & \vdots &       & \ddots& 0\cr
         0      &   0    &       &    0  & B_r
        \end{array}
  \right).
$$
A polynomial $P$ of a matrix variable 
$X \in M_{n,d}$ which is invariant under 
the action of $K$ is called $K$-invariant.
If we write $X = (X_1,\dots,X_r)$ with $X_i\in M_{n,d_i}$ ($1\le i\le r$)
then, by virtue of H. Weyl, for each $K$-invariant polynomial $P$
there exists a polynomial $Q$ such that 
$P(X) = Q(X_1 X_1^t,\dots,X_r X_r^t)$.
(Here we have used the assumption $d_i\ge n$ ($1\le i \le r$).)
Following ref.\ \cite{Ibukiyama} we call $Q$ the associated polynomial
(map) of $P$.
Note that $\GL{n}$ acts on associated polynomials $Q$ by 
mapping $Q(R_1,\dots,R_r)$ to $Q(A^t R_1 A,\dots,A^t R_r A)$ 
($A\in \GL{n}$). 

Let us introduce some notations for the certain spaces of polynomials.
Denote the space of all homogeneous polynomials $P$ of 
weight $v$ which are $K$-invariant by $\P_{n,v}(d_1,\dots,d_r)$.
Furthermore, let $\Q_{n,v}(r)$ be the space of all polynomials $Q$ in 
(the coefficients of) symmetric matrices $R_1,\dots,R_r\in M_{n,n}$ 
such that 
$Q(A^t R_1 A, \dots, A^t R_r A) = \det(A)^v Q(R_1,\dots,R_r)$
for all $A\in \GL{n}$. 
Finally, let $\Hq_{n,v}(d_1,\dots,d_r)$ be the subspace of 
$\Q_{n,v}(r)$ consisting of elements $Q$ such that 
$Q(X_1 X_1^t,\dots, X_r X_r^t)$ is pluri-harmonic for 
$X= (X_1,\dots, X_r)$ where $X_i\in M_{n,d_i}$ ($1\le i\le r$).
When $d_i\ge n$ then $\Hq_{n,v}(d_1,\dots,d_r)$ is the space of 
polynomials associated to the invariant pluri-harmonic polynomials 
in $\P_{n,v}(d_1,\dots,d_r)$.

Before we specialise to the case $r=2$ on which we will concentrate 
in the following sections let us recall the connection between 
pluri-harmonic polynomials and differential operators.

For $Z=(z_{r,s})\in\H_n$ we define $\partial_{Z}$ to be  the 
$n\times n$ matrix with components 
$$  \frac{1}{2}(1+\delta_{r,s})\frac{\partial}{\partial z_{r,s}}
$$

Let $Q$ be a polynomial in $\Q_{n}(r)$ and  
set $D = Q(\partial_{Z_1},\dots,\partial_{Z_r})$.
For fixed, even $d_{i} = 2k_{i}$ ($1 \leq i \leq r$) and fixed 
non-negative $v$, {\it i.e.} $v \in\N_0$, we say that $D$ satisfies the 
`commutation relation' if the following condition is satisfied.

\noindent
{\it For all holomorphic functions 
$F(Z_{1},\ldots,Z_{r}):\H_{n}^{r}\to\C$ 
and all $\gamma \in \SPR$ one has}
\[
D\left( F(Z_{1},\ldots,Z_{r})\prod_{i=1}^{r}\det(CZ_{i}+D)^{-k_{i}}
 \right) |_{Z=Z_{1}=\cdots=Z_{r}}
=
DF(\gamma Z,\ldots,\gamma Z) \det(CZ+D)^{-k-v}
\] 
{\it where $k= k_1+\dots+k_r$ and 
$\gamma =\left( \begin{array}{cc} A & B \\ C & D \end{array}\right)$
with $n\times n$ matrices $A,B,C$ and $D$.}

Using this condition we can now state (a special case of) Theorem 2
of ref. \cite{Ibukiyama}.
\begin{thm}\cite{Ibukiyama}
\label{Ibuthm}
For fixed integers $k_{1},\ldots,k_{r}$ with $2k_i\ge n$ 
and a non-negative integer $v$ the above commutation relation 
is satisfied if and only if 
$Q \in \Hq_{n,v}(2k_{1},\ldots,2k_{r})$. 
In particular, 
let $F_i:\H_n\to\C$ ($1\le i\le r$) be Siegel modular forms on some 
discrete subgroup $\Gamma\subset\SPR$ of finite co-volume 
with weights $2 k_i \ge n$, respectively. 
Let $Q$ be in $\Hq_{n,v}(2k_1,\dots,2k_r)$ and 
set $D = Q(\partial_{Z_1},\dots,\partial_{Z_r})$.
Then the restriction of 
$$  D\left( F_1(Z_1)\cdots F_n(Z_r) \right) $$
to $Z = Z_1=\dots=Z_r$ is a Siegel modular form 
of weight $k_1+\dots+k_r+v$ on $\Gamma$.
\end{thm}

{}From this theorem it is obvious that every element of 
$\Hq_{n,v}(d_1,\dots,d_r)$ defines a $r$-linear covariant 
differential operator and vice versa. Here a differential 
operator $D$ is called covariant if $D$ satisfies the commutation 
relation above.
In the rest of this paper we will be interested in the case 
of bilinear differential operators, {\it i.e.} the case $r=2$.
 
Finally, let us recall a result on  the dimension of 
$\Hq_{n,v}(d_1,d_2)$ which was already stated in ref.\ \cite{Ibukiyama}.
\begin{prop}\cite{Ibukiyama}
\label{Hdim}
For $d_1\ge n$ and $d_2\geq n$ one has 
$\dim(\Hq_{n,v}(d_1,d_2)) = 1$ 
for all non-negative, even $v$. 
\end{prop}
{\it Proof.} The proof for the case $d_1,d_2\ge 2n$ is already
contained in ref.\ \cite{Ibukiyama}; we give the complete proof 
in section \ref{proofs}.

\bigskip\noindent
In the next section we will, for $d_1, d_2 \geq n$, give a (quite 
explicit) description of  $\Hq_{n,v}(d_1,d_2)$. 

\section{Explicit description of $\Hq_{n,v}(d_1,d_2)$}
\label{main}

In this section we will obtain a description of 
$\Hq_{n,v}(d_1,d_2)$ for $d_1\ge n$ and $d_2 \geq n$.
As $\Hq_{n,v}(d_1,d_2)$ is one dimensional in this case 
we only have to find a non-zero element in this space.
To find such an element we study the structure of 
$\Q_{n}(2)$ first. 
More precisely, we show that $\Q_{n}(2)$ is freely generated by 
certain homogeneous polynomials $P_\alpha$ ($0\le\alpha\le n$)
of weight 2 (Lemma \ref{hompoly}). 
Using this explicit description of $\Q_{n}(2)$ we can write
an element $Q_{n,v}\in \Hq_{n,v}(d_1,d_2)$ as a linear combination 
of monomials in the $P_\alpha$'s. We then calculate the Laplacian
of $Q_{n,v}$ and find a set of recursion relations for the
coefficients in the linear combination expressing $Q_{n,v}$
in terms of the $P_\alpha$'s (Theorem \ref{vgeneral}). 
These recursion relations allow to solve uniquely 
(up to multiplication by a scalar) for the coefficients.
In the case of general $n$ and $v$, however, we do not 
obtain explicit closed formulas for the coefficients; several special 
cases where we find such formulas are discussed in \S\ref{examples}.

In order to present our line of arguments in a transparent way 
we will postpone the proofs of various lemmas and our main theorem
to the next section.

\bigskip\noindent
We can give $\Q_{n}(2)$ (and therefore 
also each $\P_{n,v}(d_1,d_2)$) explicitly as follows. 

Firstly, it is clear that $\Q_{n,v}(2) = 0$ if $v$ is odd. 
Indeed, any $Q(R,R') \in \Q_{n,v}(2)$ is determined by its 
values for $R = \id_{n}$ and diagonal matrices $R'$. 
Taking $A= (a_{i,j})$ as 
$a_{1,1}=-1$, $a_{i,i} = 1$ for $i \neq 1$ and $a_{i,j}=0$ for all 
$i \neq j$, we get $\det(A) = -1$ and $A^t R A = R$, 
$A^t R' A= R'$ for the above $R$ and $R'$. Hence, if $v$ is 
odd, then $Q(R,R') = (-1)^v Q(R,R')$ and $Q = 0$.  

Secondly, we give some typical elements of $\Q_{n,2}(2)$. 
For $0\le \alpha\le n$ denote by $P_\alpha$ the polynomial 
in the matrix variables $R\in M_{n,n}$ and $R'\in M_{n,n}$ 
defined by 
$$ \det(R + \lambda R' ) = 
     \sum_{\alpha=0}^{n} P_{\alpha}(R,R')\ \lambda^{\alpha}.
$$
Then, each $P_\alpha$ obviously belongs to $\Q_{n,2}(2)$ and, therefore, 
all polynomials in the $P_\alpha$  belong to $\Q_n(2)$. 

Thirdly, we show that the converse also holds true so that we obtain
an explicit description of $\Q_n(2)$. 
\begin{lem}
\label{hompoly}
The ring $\Q_n(2)$ is generated by the algebraically independent 
polynomials $P_0,\dots,P_n$.
\end{lem}
 
Now let $Q_{n,2v}$ be an element of 
$\Hq_{n,2v}(d_1,d_2)\subset\Q_{n,2v}(2)$.
Then, by the last lemma, $Q_{n,2v}$ can be written in the form
$$ Q_{n,2v} = \sum_{a\in I_{n,2v}} C(a)\ 
              \prod_{\alpha=0}^n P_{\alpha}^{a_\alpha}$$
where $I_{n,2v} = \{ a=(a_0,\dots,a_n)\in \N_0^{n+1}\  
                  |\ \sum_{\alpha=0}^n a_\alpha = v \}$.

The definition of $\Hq_{n,2v}(d_1,d_2)$ in the last section was rather 
indirect as $\Hq_{n,2v}(d_1,d_2)$ was defined of the space of
polynomials $Q\in \Q_{n,2v}(2)$ such that $Q(X_1 X_1^t,X_2 X_2^t)$ is 
pluri-harmonic for $X= (X_1,X_2)$ with $X_i\in M_{n,d_i}$ ($1\le i\le 2$).
One can, however, also give a direct definition as follows.

Define two differential operators $L^{(d_1)}_{i,j}$ and 
${L'}^{(d_2)}_{i,j}$ 
acting on elements of $\Q_{n,v}(2)$ by
\begin{eqnarray*}
L^{(d_1)}_{i,j} &=& d_{1}(1+\delta_{ij})D_{ij} + 4 R_{i,j} D_{i,i} D_{j,j}  +\\
         &&\qquad 
               \sum_{m' \neq j, m \neq i} R_{m',m} D_{m',j} D_{m,i}
            + 2\sum_{m' \neq j} R_{m',i} D_{m',j} D_{i,i}
            + 2\sum_{m  \neq i} R_{m ,j} D_{j ,j} D_{m,i},\\
{L'}^{(d_2)}_{i,j} & = & d_{2}(1+\delta_{ij})D'_{ij} +
                       4 R'_{ij} D'_{ii} D'_{jj} + \\
         &&\qquad 
              \sum_{m' \neq j, m \neq i} R'_{m',m} D'_{m',j} D'_{m,i}
              + 2\sum_{m'\neq j} R'_{m',i} D'_{m',j} D'_{i,i}
              + 2\sum_{m \neq i} R'_{m ,j} D'_{ j,j} D'_{m,i}
\end{eqnarray*}
where we have used $R = (R_{i,j}),R'=(R'_{i,j})\in M_{n,n}$ 
for the two matrix variables
and   $D_{i,j}$ and $D'_{i,j}$ for the differential operators
$$ D_{i,j}  = \frac{\partial}{\partial R_{ij}}, \qquad
   D'_{i,j} = \frac{\partial}{\partial R'_{ij}} \qquad (1\le i,j\le n).
$$

It is now easy to see that the operators $L^{(d_1)}_{i,j}$ and 
${L'}^{(d_2)}_{i,j}$ 
describe the action of $\Delta_{i,j}(X)$ and $\Delta_{i,j}(X')$, 
respectively on associated polynomials. More precisely, let 
$Q\in\Q_{n,v}(2)$ be the associated polynomial of a polynomial  
$P\in\P _{n,v}(d_1,d_2)$, {\it i.e.} $Q(X X^t, X' {X'}^t) = P(X,X')$
with $X\in M_{n,d_1}$ and $X'\in M_{n,d_2}$, then one has
\begin{eqnarray*}
\Delta_{i,j}(X) P(X,X')  &=& (L^{(d_1)}_{i,j} Q)(X X^t, X' {X'}^t)\\
\Delta_{i,j}(X') P(X,X') &=& ({L'}^{(d_2)}_{i,j} Q)(X X^t, X' {X'}^t).
\end{eqnarray*}

By the remarks in the last section a ``homogeneous'' polynomial is 
pluri-harmonic if and only if it is harmonic. 
Hence a necessary and sufficient condition for a polynomial  
$Q_{n,2v}\in \Q_{n,2v}(2)$ being in $\Hq_{n,2v}(d_1,d_2)$ 
is given by
$$ \sum_{i=1}^n (L^{(d_1)}_{i,i} + {L'}^{(d_2)}_{i,i})  Q_{n,2v} = 0.$$

In order to evaluate this equation we have to calculate
the action of the Laplacians $L^{(d_1)}_{i,i}$ and 
${L'}^{(d_2)}_{i,i}$ on products of $P_{\alpha}$'s.
Let us first describe the action of these operators 
on a product of two polynomials $Q$ and $Q'$
\begin{eqnarray*}
L^{(d_1)}_{i,i}(Q Q')  & = & (L^{(d_1)}_{i,i}Q) Q'  + 
                         Q (L^{(d_1)}_{i,i}Q') + 8(Q,Q')_{i,R}, \\
{L'}^{(d_2)}_{i,i}(Q Q') & = & ({L'}^{(d_2)}_{i,i}Q) Q' + 
                         Q({L'}^{(d_2)}_{i,i}Q') + 8(Q,Q')_{i,R'}
\end{eqnarray*}
where 
\begin{eqnarray*}
 4 (Q,Q')_{i,R}   &=& R_{i,i} (D_{i,i}Q) (D_{i,i} Q')+ 
                \sum_{1\le l\le n}     R_{l,i} (D_{l,i}Q) (D_{i,i} Q') +\\
      &&\qquad
               +\sum_{1\le m\le n}     R_{m,i} (D_{i,i}Q) (D_{m,i} Q')
               +\sum_{1\le l,m\le n}   R_{l,m} (D_{l,i}Q) (D_{m,i} Q')\, ,\\
 4 (Q,Q')_{i,R'}   &=& R'_{i,i} (D'_{i,i}Q) (D'_{i,i} Q')+ 
                \sum_{1\le l\le n}     R'_{l,i} (D'_{l,i}Q) (D'_{i,i} Q') +\\
      &&\qquad
               +\sum_{1\le m\le n}     R'_{m,i} (D'_{i,i}Q) (D'_{m,i} Q')
               +\sum_{1\le l,m\le n}   R'_{l,m} (D'_{l,i}Q) (D'_{m,i} Q').
\end{eqnarray*}

To state the corresponding formulas for $Q=P_\alpha$ and 
$Q' = P_\beta$ let us define polynomials $(P_{i_1,\dots,i_g;j_1,\dots,j_g})$ 
depending on two matrix variables $R,R'\in M_{n,n}$
and the variable $\lambda$ as the determinant of the matrix 
$R+\lambda R'$ with the rows $i_1,\dots,i_g$ and columns 
$j_1,\dots,j_g$ removed. Furthermore, we denote by 
$(P_{i_1,\dots,i_g;j_1,\dots,j_g})_\alpha$ the coefficient 
of $\lambda^\alpha$ in $(P_{i_1,\dots,i_g;j_1,\dots,j_g})$.

Then one has the following result.
\begin{lem}
\label{deltagrad}
With the notations above one has
\begin{eqnarray*}
L^{(d_1)}_{i,i}  P_\alpha  &=& 2(d_1+1-n+\alpha)\ (P_{i;i})_\alpha 
                          \qquad\ \ (0\le \alpha\le n-1)\, ,\\
{L'}^{(d_2)}_{i,i} P_\alpha  &=& 2(d_2+1-\alpha)\qquad\, (P_{i;i})_{\alpha-1} 
                       \qquad (1\le \alpha\le n)\, ,\\
L^{(d_1)}_{i,i}  P_n       &=& {L'}^{(d_2)}_{i,i}  P_0  =  0\, ,\\
        (P_\alpha, P_\beta)_{i,R}
   &=&   P_{\alpha} (P_{i;i})_{\beta} - 
         P_{\beta+1}(P_{i;i})_{\alpha-1}+
        (P_{\alpha-1},P_{\beta+1})_{i,R}\, ,\\
      (P_\alpha, P_\beta)_{i,R'}
   &=& P_{\beta} (P_{i;i})_{\alpha-1} - 
       P_{\alpha-1}(P_{i;i})_{\beta}+
      (P_{\alpha-1},P_{\beta+1})_{i,R'}
\end{eqnarray*}
where in the two last equations $\alpha\le \beta$ and we have set 
$P_\gamma = (P_{i;i})_\gamma =0$ for $\gamma<0$ or $\gamma>n$.
\end{lem}

Using these formulas one can now calculate 
$\sum_{i=1}^n (L^{(d_1)}_{i,i}+{L'}^{(d_2)}_{i,i})$ of any 
monomial in the $P_\alpha$'s.
To be able to extract the equations which the coefficients $C(a)$ 
have to satisfy if $Q_{n,v}$ is pluri-harmonic we need to know that
the polynomials $\sum_{i=1}^n (P_{i;i})_\beta$ are linearly independent 
over  $\Q_{n}(2)$.

\begin{lem}
\label{lindep}
For any $n \geq 1$, the polynomials $\sum_{i=1}^n (P_{i;i})_\alpha$
($0\le \alpha \le n-1$) are linearly independent over $\Q_{n}(2)$.
\end{lem}

We can now use the last three lemmas to calculate the equations 
satisfied by the coefficients $C(a)$.
\begin{thm}
\label{vgeneral}
Let $d_1\ge n$ and $d_2\ge n$. With the notations as above the polynomial  
$$ Q_{n,2v} = \sum_{a\in I_{n,2v}} C(a)\ 
   \prod_{\alpha=0}^n P_{\alpha}^{a_\alpha}
$$
is a pluri-harmonic of weight $2v$, {\it i.e.} 
$Q_{n,2v}$ is in $\Hq_{n,2v}(d_1,d_2)$,
if the coefficients $C(a)$ satisfy
\begin{eqnarray*}
  (d_1+1-n+i+2(a_i-1)) a_i\        C(a)             &&= 
 -(d_2-i)             (a_{i+1}+1)\ C(a-e_i+e_{i+1}) \\
 && -2\sum_{{i<l\le l'}\atop{l'+l-i\le g}}
     \tilde a(i,l,l')_l\ ( \tilde a(i,l,l')_{l'}-\delta_{l,l'})\ 
     C( \tilde a(i,l,l') )\\
 && +2\sum_{{i<l\le l'}\atop{l'+l-i-1\le g}}
     \hat a(i,l,l')_l\ (\hat a(i,l,l')_{l'}- \delta_{l,l'})\
     C( \hat a(i,l,l') ) 
\end{eqnarray*}
for  all $a\in I_{n,2v}$ such that $i:=\min\{j | a_j \not=0\} < n$. 
Here we have used $\tilde a(i,l,l') = a-e_i+e_l+e_{l'}-e_{l+l'-i}$,
$\hat a(i,l,l') = a-e_i+e_l+e_{l'}-e_{l+l'-i-1}$ and 
$e_j\in \{0,1\}^{n+1}$ for the vector with components 
$(e_j)_l =  \delta_{j,l}$. Note that in the above formula we have 
set $C(a)=0$ if $a_j<0$ for some $0\le j\le n$.
\end{thm}

Firstly, note that the equations determine the coefficients 
$C(a)$ uniquely as multiples of $C((0,\dots,0,v))$.
To see this define an order for the elements of $I_{n,v}$ 
by saying that $a < b$ for some $a,b\in I_{n,v}$ if there exists 
some $j$ with $0\le j\le n$ so that $a_j<b_j$ and $a_i=b_i$ for all $0\le i<j$.  
Looking at the above recursion equations it is obvious that, for $d_1\ge n$, 
they can be used to express $C(a)$ in terms of coefficients $C(b)$ with
$b<a$. Hence we obtain inductively that once given $C((0,\dots,0,v))$
all other coefficients are uniquely determined (assuming still $d_1\ge
n$). This implies that there is exactly one solution to the equations
which is indeed easy to calculate for any given {\it numerical} values
of $n$ and $v$. It would of course be desirable to give closed 
explicit formulas for the coefficients $C(a)$ in general but we have 
not succeeded in doing so; several special cases where we obtain
such formulas are discussed in \S\ref{examples}.

As, by Proposition \ref{Hdim}, the dimension of $\Hq_{n,2v}(d_1,d_2)$ 
is equal to one for $d_1,d_2\ge n$  Theorem \ref{vgeneral} gives a 
(more or less explicit) description of $\Hq_{n,2v}(d_1,d_2)$.
\medskip

Finally, the differential operators obtained from non-constant
invariant pluri-harmonic polynomials in $\Hq_{n,v}(d_1,d_2)$ 
map two modular forms to a cusp form.
\begin{prop}
\label{cusp}
Let $\Gamma$ be a subgroup of $\SPQ$ which is commensurable with 
$\SP$ and let $F$ and $G$ be Siegel modular forms on $\Gamma$ of weight
$k$ and $l$, respectively.
Let $D$ be a covariant differential operator obtained from a 
non-constant invariant pluri-harmonic polynomial in $\Hq_{n,v}(2k,2l)$. 
Then $D(F(Z_{1})G(Z_{2}))|_{Z_{1}=Z_{2}}$ is a cusp form.
\end{prop}

\medskip

In section \ref{proofs} we will give the proofs of the Lemmas
\ref{hompoly}, \ref{deltagrad} and \ref{lindep} as well as 
the proof of Theorem \ref{vgeneral} and  Proposition \ref{cusp}.

\section{Proofs}
\label{proofs}
In this section we have collected those proofs which have not been 
given so far. 

\subsection{Proof of Proposition \ref{Hdim}}
\label{Hdimproof}

{\it Proof of Proposition \ref{Hdim}.}
Roughly speaking, the irreducible representations
of $O(d)$ are parametrised by ``Young diagrams'' 
$(f_{1},...,f_{k})_{+}$ and 
$(f_{1},...,f_{k})_{-}$ where
$f_{1} \geq f_{2} \geq \cdots \geq f_{k} \geq 0$,
$k = d/2$, and, where,  $+$ and $-$ coincide if $f_{k} \neq 0$.
The space of pluri-harmonic polynomials $P(X,X^{'})$ 
($(X,X^{'}) \in M_{n,d_{1}} \times M_{n,d_{2}}$)
such that $P(A^t X,A^t X^{'}) = \det(A)^{v}P(X,X^{'})$ gives an 
irreducible representation of $O(d)$ ($d = d_{1} + d_{2}$)  
corresponding to the Young diagram $(v,...,v)_{+}$ 
of depth $n$ ({\it cf.} the notation of M.\ Kashiwara and M.\ Vergne
in ref.\ \cite{KashiwaraVergne}). 
If we take the restriction of this representation to 
$O(d_{1}) \times O(d_{2})$ then $\Hq_{n,v}(d_{1},d_{2})$ is the 
subspace which corresponds to the trivial representation of 
$O(d_{1}) \times O(d_{2})$. So, what we should do is to count the 
multiplicity of the trivial representation in the restriction of 
$(v,...,v)_{+}$ to  $O(d_{1}) \times O(d_{2})$.
The irreducible decomposition of a similar restriction  
has already been worked out by K.\ Koike and I.\ Terada in
Theorem 2.5 and Corollary 2.6 on p.\ 115 of ref.\ \cite{KoikeTerada} 
but here a subtle point is different. 
For the sake of simplicity let us denote by $R(O(d))$ the ring of
those characters of $SO(d)$ which can be obtained as restrictions of 
characters of irreducible representations of $O(d)$. 
K.\ Koike and I.\ Terada  take $\lambda_{SO(d)} \in R(O(d))$
and give the following formula for the restriction of
$\lambda_{SO(d)}$ to $O(d_{1}) \times O(d_{2})$.
\[
\lambda_{SO(d)} 
= \sum_{\beta, \mu, \kappa, \nu}
L^{\lambda}_{\beta,\mu}L^{\beta}_{2\kappa,\nu}
\pi_{SO(d_{1})}(\mu_{SO})\times \pi_{SO(d_{2})}(\nu_{SO})),
\]
where $L^{a}_{b,c}$ are the so-called Littlewood-Richardson coefficients, 
which can be calculated explicitly in principle, where
the parameters $\beta$, $\kappa$, $\mu$, $\nu$ are any partitions or
 ``universal characters'' and, where  $\pi_{SO(c)}$ is a ``specialisation'' 
homomorphism  whose image is contained in $R(O(c))$ ({\it cf. loc. cit.}).
Then what we should do is as follows. 
First we calculate the coefficients 
for those pairs $\mu$ and $\nu$ where 
$\pi_{SO(d_{1})}(\mu)$ and $\pi_{SO(d_{1})}(\nu)$ are the trivial
characters.

If $d_{1} > n$ and $d_{2} > n$ it is not difficult to show 
that this occurs only when $\mu$ and $\nu$ are trivial and that  
the coefficient is one. This is seen as follows.
There are exactly two irreducible representations of $O(c)$ whose 
restriction to $SO(c)$ is trivial:
the trivial representation and the determinant representation. 
We must exclude the latter possibility. Since we are assuming 
$d_{1} > n$ and $d_{2} > n$ the fundamental theorem on invariants 
({\it cf.} Theorem 2.9A on p.\ 53 of ref.\ \cite{Weyl}) implies that 
any $SO(d_{i})$ invariant vector is also $O(d_{i})$ invariant. 
This means in our case that if the restriction to $SO(d_{i})$ is 
trivial then it comes from the trivial representation of $O(d_{i})$.    
This proves our assertion for $d_{1} > n$ and $d_{2} > n$.
 
When $d_{1} = n$ or $d_{2} = n$ the proof is more involved. 
First of all, for $\lambda = (v,\ldots,v)_{+}$,  
we pick up pairs of $\mu$ and $\nu$ such that 
$L^{\lambda}_{\beta,\mu}L^{\beta}_{2\kappa,\nu} \neq 0$ 
and $\pi_{SO(d_{1})}(\mu) \times \pi_{SO(d_{2})}(\nu)$ is 
trivial. Under the assumption that $d_{1} \geq n$ and $d_{2} \geq n$ 
these pairs can be described as follows
(for the sake of simplicity we set $\rho_{a}=(a,\ldots,a)$ 
 with depth $n$). 
\\
(1) If $v$ is even and $\mu$ and $\nu$ are trivial then the 
    coefficient is one.
\\
(2) If $d_{1} = n$ and $v$ is odd then the coefficient is one for 
    $\mu = \rho_{1}$, $\nu$ trivial, $\beta=\rho_{v-1}$ and   
    $\kappa = \rho_{(v-1)/2}$.
\\
(3) If $d_{2} = n$ and $v$ is odd then the coefficient is one for 
    $\mu$ trivial, $\nu = \rho_{1}$, $\beta = \lambda$ and
    $\kappa=\rho_{(v-1)/2}$.
\\
(4) If $d_{1} = d_{2} = n$ and $v$ is even with $v \geq 2$ 
    then the coefficient is one for $\mu=\rho_{1}$, $\nu=\rho_{1}$, 
    $\beta=\rho_{v-1}$ and $\kappa=\rho_{(v-2)/2}$.
\\
These possibilities exhaust all cases giving the trivial
representation of $SO(d_{1}) \times SO(d_{2})$. 
We have already shown in \S\ref{main} that $v$ is even in our case 
so that the cases (2) and (3) do not occur. 

If $d_{1} > n$ or $d_{2} > n$ then the pluri-harmonic polynomial 
in question is invariant by $O(d_{1})$ or $O(d_{2})$, respectively. 
Indeed, if this were not the case then, by the theory of invariants, it is of the form $\det(X^{'})Q(R,R^{'})$ ($d_{2} = n$ and $X^{'} \in M_{n}$)
or $\det(X)Q(R,R^{'})$ ($d_{1} = n$, $X \in M_{n}$), respectively, where 
$Q \in \Q_{n,v-1}(2)$. But since $Q = 0$ unless $v-1$ is even this 
cannot be the case. 

Therefore, the only remaining case is $d_{1} = d_{2} = n$ and we must show
that the trivial representation of $O(n) \times O(n)$ occurs exactly
once. We prove this by showing that there exists exactly one 
$\det(g) \times \det(h)$ representation of $O(d_{1}) \times O(d_{2})$ 
in our restriction. Since, for $d_{1} = d_{2} = n$ and $v$ even, the
multiplicity of the trivial representation of 
$SO(d_{1}) \times SO(d_{2})$ is two this indeed implies that our
restriction contains the trivial representation of 
$O(n) \times O(n)$ exactly once.

For even $v$ assume that $P(X,X^{'})$ ($X$, $X^{'} \in M_{n}$) is a 
pluri-harmonic polynomial such that 
$P(A^t X,A^t X^{'}) = \det(A)^{v}P(X,X^{'})$. 
By the above considerations this polynomial is invariant both by 
$O(d_{1})$ and $O(d_{2})$ or odd invariant for both $O(d_{1})$ and 
$O(d_{2})$, {\it i.e.} $P(Xg,X^{'}h) = P(X,X^{'})$ or
$P(Xg,X^{'}h) = \det(g)\det(h)P(X,X^{'})$, respectively.
We show the latter case occurs just for one polynomial in question 
(up to multiplication by a constant).
By the classical theorem of invariants (Weyl, {\it loc.cit.})
we find in the latter case $P(X,X^{'}) = \det(X)\det(X^{'})Q(R,R^{'})$ 
for some $Q \in \Q_{n,v-2}(2)$. By applying 
$\Delta_{11} = \Delta_{11}(X) + \Delta_{11}(X^{'})$ 
to this polynomial it is easy to see that 
\[
\Delta_{11}(\det(X)\det(X^{'})Q(R,R^{'})) 
=
(\det(X)\det(X^{'}))
(\Delta_{11}Q + 
2(\frac{\partial Q}{\partial R_{11}}
+\frac{\partial Q}{\partial R^{'}_{11}})).
\]
Hence the polynomial is pluri-harmonic if and only if 
$\Delta_{11}Q + 2(\frac{\partial Q}{\partial R_{11}}
+\frac{\partial Q}{\partial R^{'}_{11}}) = 0$.
The latter action, however, is nothing but the action of 
$L^{(d_{1}+1)}_{1,1} + L^{'(d_{2}+1)}_{1,1}$. 
Since $d_{1} + 1 > n$ and $d_{2} + 1 > n$ we have already 
shown that  the kernel of $L^{(d_{1}+1)}_{1,1} + L^{'(d_{2}+1)}_{1,1}$
is one dimensional. 
\qed

\subsection{Proof of Lemmas \ref{hompoly}-\ref{lindep}, 
            Theorem \ref{vgeneral} and 
            Proposition \ref{cusp}}

This subsection contains the proofs of the three lemmas 
stated in section \ref{main} which we then use to prove Theorem 
\ref{vgeneral}. We also give a proof of Proposition \ref{cusp}.

\smallskip\noindent
{\it Proof of Lemma \ref{hompoly}.}
For any polynomial $Q(R,R') \in \Q_{n,v}(2)$
the polynomial $Q(\id_{n},T)$ is a polynomial of the 
coefficients of $T$ where $T$ is symmetric and 
$\id_{n}$ is the unit matrix of size $n$. 
Since $Q(\id_{n},T) = Q(\id_{n},O^{-1} T O)$ for any orthogonal 
matrix $O$ the polynomial $Q(\id_{n},T)$ is a polynomial in 
the functions $\mu_{i}(T)$ ($1 \leq i \leq n$) defined by
$\det(t \id_n + T) = \sum_{i=0}^n \mu_i(T) t^i$.
Indeed, if $T$ is a diagonal matrix then $Q(\id_{n},T)$
is a symmetric function of the diagonal components and hence 
a function of the $\mu_{i}(T)$. 
Since $\mu_{i}(T)$ is invariant under conjugation of $T$ we obtain 
that $Q(\id_{n},T)$ is a polynomial in the $\mu_{i}(T)$.
Furthermore, for $0 \leq i \leq n$, we find that 
$\det(R)\mu_i(T) = P_i(R,R')$ 
where we have set $R= A^t A$ and $T = (A^t)^{-1}R' A^{-1}$
(strictly speaking we are working here over an algebraic extension 
of our ring $\C[R_{i,j}]$ ($R = (R_{i,j})$) 
which allows to find such a matrix $A$ with $R=A^t A$).
Hence, for any non-negative integers $e_{i}$ ($1\leq i \leq n$) and 
$l_{0} = \sum_{i=1}^{n}e_{i}$ we get 
$$\det(R)^{l_{0}}\prod_{i=1}^{n}\mu_{i}(T)^{e_{i}}
= \prod_{i=1}^{n}P_{i}^{e_{i}}(R,R'). 
$$
Now, any linear combination of the above 
``monomials'' for fixed $l_{0}$ and various $e_{i}$ ($i \geq 1$)
is not divisible by $\det(R)$ since even if some column of 
$R$ is a zero vector such a linear combination does not vanish.
This is seen as follows. Choose $R' = \id_{n}$
and $R$ as a diagonal matrix with $R_{1,1} = 0$. 
Then, each $P_{i}$ ($ 1 \leq i \leq n-1$) becomes an elementary symmetric 
function of the $R_{i,i}$ ($2 \leq i \leq n$). 
This means that the above linear combination does not vanish identically.
Furthermore, it also means that, for a fixed $v$, any linear combination 
of the ``monomials'' $\det(R)^{l}\prod_{i=1}^{n}P_{i}^{e_{i}}(R,R')$
with $l + \sum_{i=1}^{n}e_{i} = v$ where $l$ is an integer and 
the $e_i$ are non-negative integers is a polynomial if and only if 
$l$ is positive.
Since any element in $Q_{n,v}(2)$ is a linear combination 
of ``monomials'' of the above type $Q_{n,v}(2)$ is generated by the 
$P_\alpha$ ($0\le \alpha\le n$). Furthermore, 
the $P_i$  ($ 1 \leq i \leq n$) are obviously algebraic independent
with each other. This completes the proof of the lemma.\qed

\bigskip\noindent
{\it Proof of Lemma \ref{deltagrad}.}
To prove the  equations in Lemma \ref{deltagrad} 
we choose without loss of generality  $i=1$.

Firstly, note that one has 
$$D_{l,1}P_{\alpha} = (2-\delta_{1,l})(-1)^{l+1}(P_{1;l})_{\alpha}
  \qquad (0 \leq \alpha \leq  n)
$$
(here we regard $(P_{1;j})_{n} = 0$) and 
$$ D_{l',1} (P_{1;l})_{\alpha} = (-1)^{l'}(P_{1,l';1,l})_{\alpha}
   \qquad (l'\not=1).
$$
This directly implies 
$$ L^{(d_1)}_{1,1} P_\alpha = 
   2 d_1 (P_{1;1})_\alpha -
   2 \sum_{l,l'\not=1} (-1)^{l+l'} R_{l,l'}\, (P_{1,l';1,l})_{\alpha}.
$$
Secondly, we show that 
$$\sum_{l,l'\not=1} (-1)^{l+l'} R_{l,l'}\, (P_{1,l';1,l})_{\alpha}
 = (n-1-\alpha) (P_{1;1})_\alpha.
$$
This obviously implies the formula stated in the lemma.

The last equality can be proven as follows.
Multiplying both sides with $\lambda^\alpha$ and summing over $\alpha$
gives, as an equivalent equation, 
$$ \lambda \frac{d}{d\lambda} P_{1;1} =
    (n-1)\ P_{1;1} - 
    \sum_{j,k\not=1} (-1)^{j+k}R_{j,k} P_{1,j;1,k}.
$$
Instead of proving this equation is suffices to prove the equation with 
$P_{1;1}$ replaced by $P$ and $n-1$ replaced by $n$, 
{\it i.e.} to prove
$$ \lambda \frac{d}{d\lambda}  P = n\,  P - 
    \sum_{j,k=1}^n (-1)^{j+k}\ R_{j,k}  P_{j;k}.
$$
Expanding the determinants on the right hand side gives
\begin{eqnarray*} 
&& n\sum_{\sigma\in S_n} (-1)^\sigma \prod_{i=1}^{n}
        (R_{i,\sigma(i)} + \lambda R'_{i,\sigma(i)})\\ 
&& -  \sum_{k=1}^n \sum_{\sigma\in S_n} 
      (-1)^\sigma 
      \prod_{i=1}^{n}
        (R_{i,\sigma(i)} + \lambda (1-\delta_{i,k}) R'_{i,\sigma(i)}) \\
&=&        
\sum_{k=1}^n \sum_{\sigma\in S_n} 
             (-1)^\sigma
             \lambda R'_{k,\sigma(k)}
             \prod_{i\not=k}
             (R_{i,\sigma(i)} + \lambda R'_{i,\sigma(i)})
\end{eqnarray*}
where $S_n$ is the symmetric group of $n$ elements.
Note that the last expression on the r.h.s. is equal to 
$\lambda \frac{d}{d\lambda} P$ so that we have proven the 
desired equality.

The analogous equation for ${L'}^{(d_2)}_{i,i}$ follows directly from 
the symmetry $R \leftrightarrow R'$, $\alpha \leftrightarrow n-\alpha$.

\bigskip\noindent
We now prove the equation equation for $(P_\alpha,P_\beta)_{1,R}$.
Using again the equality  
$D_{l,1}P_{\alpha} = (2-\delta_{1,l})(-1)^{l+1}(P_{1;l})_{\alpha}$
we immediately get 
$$ (P_{\alpha},P_{\beta})_{1,R} 
    = 
    \sum_{l,m=1}^{n}(-1)^{l+m}R_{l,m}(P_{1;l})_{\alpha}(P_{1;m})_{\beta}.
$$
Furthermore, we find  
$$ \sum_{m=1}^{n}(-1)^{l+m}R_{l,m}(P_{1;m})_{\alpha} 
    = 
    \delta_{1,l}P_{\alpha}- \sum_{m=1}^{n}(-1)^{l+m}
                               R_{l,m}^{'}(P_{1;m})_{\alpha-1}
$$
and 
$$ \sum_{l=1}^{n}(-1)^{l+m}R_{l,m}^{'}(P_{1;l})_{\beta}
   = \delta_{m,1}P_{\beta+1} -
      \sum_{l=1}^{n}(-1)^{l+m}R_{l,m}(P_{1;l})_{\beta+1}.
$$
Finally, collecting the above formulas gives the desired 
equation
\begin{eqnarray*}
(P_{\alpha},P_{\beta})_{1,R}
 & = & 
 \sum_{l=1}^{n}\delta_{l,1}P_{\alpha}(P_{1;l})_{\beta}
-\sum_{l,m=1}^{n}(-1)^{l+m}R^{'}_{l,m}(P_{1;m})_{\alpha-1}(P_{1;l})_{\beta}
 \\
 & = & 
P_{\alpha}(P_{1;1})_{\beta}
\\
&& \qquad -
\sum_{m=1}^{n}\left(\delta_{m,1}P_{\beta+1}(P_{1;m})_{\alpha-1}
-\sum_{l=1}^{n}(-1)^{l+m}R_{l,m}(P_{1;m})_{\alpha-1}(P_{1;l})_{\beta+1}\right)
\\ 
 & = & 
P_{\alpha}(P_{1;1})_{\beta} - P_{\beta+1}(P_{1;1})_{\alpha-1}
+ \sum_{l,m=1}^{n}(-1)^{l+m}R_{l,m}(P_{1;m})_{\alpha-1}(P_{1;l})_{\beta+1}
\\
 & = & 
P_{\alpha}(P_{1;1})_{\beta} - P_{\beta+1}(P_{1;1})_{\alpha-1}
+ (P_{\alpha-1},P_{\beta+1}).
\end{eqnarray*}

The corresponding equation for $(P_\alpha, P_\beta)_{1,R'}$  follows,
again, directly from the symmetry $R \leftrightarrow R'$, 
$\alpha\leftrightarrow n-\alpha$. 
\qed

\bigskip\noindent
{\it Proof of Lemma \ref{lindep}.}
For the sake of simplicity we denote by $(P_{i;i})_{\beta}^{n}$ 
the differentiation of $(P_{i;i})_{\beta}$ by
$\frac{\partial }{\partial R_{nn}}$. (In other words, 
$(P_{i;i})_{\beta}^{n}$ is obtained from $(P_{i;i})_{\beta}$ by 
substituting $R_{n,j} = 0$ ($n \neq j$), $R_{n,n}=1$ and
$R^{'}_{n,j} = 0$ for all $j$ with $1 \leq j \leq n$.)

For $\beta$ with $0 \leq \beta \leq n-1$, 
we take polynomials $Q_{\beta}$ of $n+1$ variables and assume that 
they satisfy
$$
\sum_{\beta=0}^{n-1}Q_{\beta}(P_{0},...,P_{n})
    (\sum_{i=1}^{n}(P_{i;i})_{\beta})
 = 0.
$$
We show that $Q_{\beta} = 0$ for all $\beta$. 
We set the last row (and the last column) of $R^{'}$ to 0 and 
also set $R_{n,j} = R_{j,n} = 0$ for all $j \neq n$. 
Then $P_{n}$ becomes $0$ and $(P_{n;n})_{\beta}$ is unchanged for 
any $\beta$. Furthermore, $P_{\alpha}$ becomes equal to 
$R_{n,n} (P_{n;n})_{\alpha}$ for $\alpha\le n-1$ 
and $(P_{i;i})_{\alpha}$ becomes equal to 
$R_{nn} (P_{i,i})_{\alpha}^{n}$  for $i\le n-1$. 
Hence we have
$$ \sum_{\beta=0}^{n-1}
     Q_{\beta}(R_{n,n} (P_{n;n})_{0},...,R_{n,n}(P_{n;n})_{n-1},0)
    (\sum_{i=1}^{n-1} R_{n,n} (P_{i;i})_{\beta}^{n} + 
    (P_{n;n})_{\beta}) = 0.
$$
Note that $(P_{n;n})_{\beta}$ does not contain $R_{n,n}$ so that 
by comparing the degree of $R_{n,n}$ on both sides we obtain  
$$ \sum_{\beta=0}^{n-1}
   Q_{\beta}(R_{n,n} (P_{n;n})_{0},...,R_{n,n} (P_{n;n})_{n-1},0)
   (\sum_{i=1}^{n-1} R_{n,n} (P_{i;i})_{\beta}^n)  = 0.
$$
Now we set $R_{n,n}=1$ and use induction by $n$. 
Since for $n=2$ the lemma is easily verified we can assume that 
the lemma holds for $n-1$. Then the above induction step gives  
nothing but the relation for $n-1$ and  we obtain 
$Q_{\beta}(X_{0},...,X_{n-1},0) = 0.$
Therefore, all $Q_{\beta}$ are divisible by $X_{n}$. 
By repeating this process we find that $Q_{\beta} = 0$.
\qed

\bigskip\noindent
{\it Proof of Theorem \ref{vgeneral}.}
Using Lemma \ref{deltagrad} we can calculate the Laplacian of $Q_{n,v}$ and
obtain, since the $\sum_{i=1}^n (P_{i;i})_\alpha$ ($0\le \alpha \le n-1$) 
are linearly independent over $\Q_{n}(2)$ by Lemma \ref{lindep},
the following set of equations for the coefficients $C(a)$
\begin{eqnarray*}
(d_1+1-g+i)a_i C(a) &+& (d_2-i)(a_{i+1}+1) C(a-e_i+e_{i+1}) =\\
&& -2\sum_{l\le l'\le i} \tilde a(i,l,l')_l\ 
                      ( \tilde a(i,l,l')_{l'}-\delta_{l,l'})\ 
                       C( \tilde a(i,l,l') )\\
&& -2\sum_{i\le l\le l'} \hat a(i,l,l')_l\ 
                      ( \hat a(i,l,l')_{l'}-\delta_{l,l'})\ 
                       C( \hat a(i,l,l') )\\
&& +2\sum_{i<l\le l'   } \tilde a(i,l,l')_l\ 
                      ( \tilde a(i,l,l')_{l'}-\delta_{l,l'})\ 
                       C( \tilde a(i,l,l') )\\
&&+2\sum_{l\le l'<  i} \hat a(i,l,l')_l\ 
                      ( \hat a(i,l,l')_{l'}-\delta_{l,l'})\ 
                       C( \hat a(i,l,l') )
\end{eqnarray*}
for all $0\le i\le n-1$.
Choosing $i=\min\{j|a_j\not=0\}$ gives the equations in the
formulation of the theorem. (These equations are equivalent to  
the vanishing of coefficient in front of 
$(\sum_{k=1}^{n} (P_{k;k})_i) \prod_{j=0,n} P_j^{a_j-\delta_{i,j}}$.)

Finally, note that by Lemma \ref{hompoly},  
Proposition \ref{Hdim} and the remarks after 
Theorem \ref{vgeneral} we know that the other 
equations, {\it i.e.} those with 
$i\not=\min\{j|a_j\not=0\}$, do not contain any 
further information about the coefficients $C(a)$.
\qed

\bigskip\noindent
{\it Proof of Proposition \ref{cusp}.}
Firstly, we show that if $R$ and $R^{'}$ are positive semi-definite 
symmetric real matrices and if $\det(R+R^{'}) = 0$ 
then $\det(R+\lambda R^{'})=0$ and   
hence $P_{\alpha}(R,R^{'}) = 0$ for all $\alpha$ with 
$0 \leq \alpha \leq n$. 
In general, if $A$ and $B$ are positive semi-definite symmetric 
matrices so is $A+B$ and $\det(A+B) \geq \det(A) \geq 0$. 
Now, let $\lambda$ be any real number between 0 and 1. Then 
we obtain $0 = \det(R+R^{'}) \geq \det(R+\lambda R^{'}) \geq 0$
so that $\det(R+\lambda R^{'})$ vanishes identically as
$\det(R+\lambda R^{'})$ is a polynomial in $\lambda$.
 
Secondly, let $F$ and $G$ be Siegel modular forms of weight $k$ and 
$l$, respectively. 
We consider $(F|^{k}_{M})(Z_{1})$ and $(G|^{l}_{M})(Z_{2})$ 
for $M \in \SPQ$. As these are modular forms on $M^{-1}\Gamma M$ 
there exists some natural number $N$ such that they have a 
Fourier expansion of the form 
\begin{eqnarray*}
(F|^{k}_{M})(Z_{1}) & = & 
\sum_{T_{1}}a_{1}(T_{1}) \exp(2\pi i\ \tr(T_{1}Z_{1}/N)), \\
(G|^{l}_{M})(Z_{2}) & = & 
\sum_{T_{2}}a_{2}(T_{2}) \exp(2\pi i\ \tr(T_{2}Z_{2}/N)),
\end{eqnarray*}
where $T_{1}$ and $T_{2}$ run over all positive semi-definite 
half-integral matrices. 
We set $H = (F|^{k}_{M})(Z_{1})(G|^{l}_{M})(Z_{2})$.
Assume now that $D$ is obtained from an associated polynomial 
$Q\in \Hq_{n,v}(2k,2l)$. 
Then we easily obtain  
\[
(DH)|_{Z_{1}=Z_{2}}=
\sum_{T_{1},T_{2}}a_{1}(T_{1})a_{2}(T_{2})
Q(T_{1}/N,T_{2}/N)\exp(2\pi i\ \tr((T_{1}+T_{2})Z_{2}/N)),
\]
where, again,  $T_{1}$ and $T_{2}$ run over all positive semi-definite
half-integral matrices, and 
$Q(T_{1}/N,T_{2}/N) = 0$ if $\det(T_{1}+T_{2}) = 0$
by the discussion above. 
By Theorem \ref{Ibuthm} we also know that
\[
(D(FG)|_{Z_{1}=Z_{2}=Z})|^{k+l+v}_{M}
 = D(H)|_{Z_{1}=Z_{2}=Z}.
\]
This means that 
$\Phi((D(FG)|_{Z_{1}=Z_{2}=Z})|^{k+l+v}_{M}) = 0$ for all $M\in\SPQ$ 
where $\Phi$ is the Siegel operator and hence  
$D(FG)|_{Z_{1}=Z_{2}=Z}$ is a cusp form. 
\qed

\section{Some explicit examples}
\label{examples}

In this section we discuss some special cases of Theorem
\ref{vgeneral}. In particular we recover the results of 
ref.\ \cite{Cohen} for $n=1$ (see also the examples in ref.\ \cite{Ibukiyama})
and  ref.\ \cite{ChoieEholzer} for $n=2$. 
Furthermore, we give closed explicit formulas for $Q_{n,v}$ for
$v=2,4$ and general $n$.

\subsection{The case $v=2$}

By the discussion in section \ref{main} we can write 
$Q_{n,2}$ as
$$ Q_{n,2} = \sum_{\alpha=0}^n C(\alpha) P_\alpha.$$
The recursion equations given in Theorem \ref{vgeneral}
simplify in this case to 
$$ (d_2-\alpha)C(\alpha) = -(d_1+n-\alpha+1) C(\alpha+1)
   \qquad (0\le \alpha \le n-1).
$$
Hence we find that, up to a constant multiple, 
$Q_{n,2}$ is given by
$$ Q_{n,2} = \sum_{\alpha+\beta=n} (-1)^\alpha\, \alpha!\, \beta!\,
                                   \bin{d_2-\alpha}{\beta}\, 
                                   \bin{d_1-\beta}{\alpha}\,
                                   P_\alpha.
$$
 
\subsection{The case $v=4$}
\label{v4}
We consider the case $v=4$ for general $n$. Writing
$$ Q_{n,4} = \sum_{i,j=0}^n C_{i,j} P_i P_j $$
with $C_{i,j} = C_{j,i}$ one finds from 
Theorem \ref{vgeneral} the following recursion equations
$$ (d_1+1-n+i + 2\delta_{i,j} ) C_{i,j} + 
   (d_2-i   + 2\delta_{i<j}      ) C_{i+1,j} =
  2\sum_{r=1}^{j-i-1} \left( C_{i+r,j-r} - C_{i+r,j+1-r} \right)
$$
where $0\le i\le j \le n$ and $i\le n-1$.

One can give an explicit solution to these equations
\footnote{This explicit form of the solution is due to D.\ Zagier}. 
The coefficients  $C_{r,s}$ ($0\le r,\,s\le n$) are given by the 
closed formula 
$$ C_{r,s}=(-1)^{r-s}\,\frac{(d_1-n+r)!(d_1-n+s)!(d_2-r)!(d_2-s)!}
     {(d_1-n)!(d_1-n+2)!(d_2-n)!(d_2-n+2)!}\,
   p_{r-s}\bigl(\kappa_1\bigl(\frac{r+s}2\bigr),\,
     \kappa_2\bigl(\frac{r+s}2\bigr)\bigr)\,,$$
where  $\kappa_1(r) = d_1+2-(n-r)$, $\kappa_2(r)= d_2+2-r$ and 
$p_e(x,y)$ is the polynomial of degree 4 in $e^2$, $x$ and $y$ given by
   $$p_e(x,y) = x^2y^2+(e^2-1)\,xy\bigl(x+y+\frac{e^2-6}6\bigr)
    +\frac{e^2(e^2-4)}{12}\,\bigl(x^2+y^2-\frac14\bigr)\,.$$
The first few values of this polynomial are given by
\begin{eqnarray*} 
   p_0(x,y)&=&x(x-1)y(y-1)\,,\\
   p_1(x,y)&=&(x+\tfrac12)(x-\tfrac12)(y+\tfrac12)(y-\tfrac12)\,,\\
   p_2(x,y)&=&xy(xy+3x+3y-1)\,,\\
   p_3(x,y)&=&(x+\tfrac12)(y+\tfrac12)
     (xy+\tfrac{15}2x+\tfrac{15}2y-\tfrac{15}4)\,,\\
   p_4(x,y)&=&x^2y^2+15x^2y+15xy^2+16x^2+25xy+16y^2-4\,.
\end{eqnarray*}
To check the correctness of the formula, it is convenient to replace
the recursion relation above by the simpler 4-term 
recurrence
  $$\kappa_1(r+3)\,C_{r,s}+\kappa_2(r+2)\,C_{r+1,s}
    =\kappa_1(r-2)\,C_{r-1,s+1}+\kappa_2(r-3)\,C_{r,s+1}\,.$$
It is easy to verify that the solution given above indeed  satisfies
this recursion relation.  

Note that the denominator in the formula for $C_{r,s}$ is of course just 
conventional. One could also make the choice $(d_1-n)!^2(d_2-n)!^2$, 
which gives the simpler expression 
  $$ C_{r,s}=(-1)^{r-s}(n-1-d_1)_r(n-1-d_1)_s(n-1-d_2)_{n-r}(n-1-d_2)_{n-s}
   p_{r-s}\bigl(\kappa_1\bigl(\frac{r+s}2\bigr),
     \kappa_2\bigl(\frac{r+s}2\bigr)\bigr)$$
(here we have used $(x)_n$ for $\prod_{i=0}^{n-1} (x-i)$).

\subsection{The genus one case}
\label{genus1}

In this case $\Q_{1}(2)$ is generated by $P_0$ and $P_1$.
Writing $C_{r,s}$ for $C(a)$ with $a=(r,s)$ the 
recursion equations in Theorem \ref{vgeneral} for $C(a)$ become
$$ 2(d_1/2+r-1)r \ C_{r,s} = - d_2 (s+1)\ C_{r-1,s+1} +
                           2 s (s+1)\ C_{r-1,s+1}
                       = -2 (d_2/2-s)(s+1)\  C_{r-1,s+1}
$$
for $1\le r \le v$.
Hence we find that, up to a constant multiple, 
$Q_{1,2v}$ is given by
$$ Q_{1,2v} = \sum_{r+s=v} (-1)^r\, \bin{v+d_2/2-1}{r}\, \bin{v+d_1/2-1}{s}\ 
                           P_0^r\, P_1^s. 
$$
The differential operators corresponding to $Q_{1,2v}$ are 
precisely the operators studied by H.\ Cohen in ref.\ \cite{Cohen}.
Note that this case has already been discussed in the 
context of pluri-harmonic polynomials in ref.\ \cite{Ibukiyama}.

\subsection{The genus two case}
\label{genus2}

To make contact with the results of ref.\ \cite{ChoieEholzer}
we use $P^*_0=P_0,P^*_2=P_2$ and $P^*_1=P_0+P_1+P_2$ 
instead of $P_0, P_1$ and $P_2$ to generate $\Q_{2}(2)$.
Writing $Q_{2,2v}$ as
$$ Q_{2,2v} = \sum_{r+s+p=v} C_{r,s,p}\, (P^*_0)^r\, (P^*_2)^s\, (P^*_1)^p
$$
and applying Lemma \ref{deltagrad} and Lemma \ref{lindep}
we find the following two recursion relations for the coefficients 
$C_{r,s,p}$ ({\it cf.} the equations on page 13 of {\it loc. cit.}) 
\begin{eqnarray*}
0 &=& (r+1)((d_1-3)/2      +r+1) C_{r+1,s,p} + 
      (p+1)((d_1+d_2-3)/2  +p+1) C_{r,s,p+1}\,, \\
0 &=& (s+1)((d_2-3)/2      +s+1) C_{r,s+1,p} + 
      (p+1)((d_1+d_2-3)/2  +p+1) C_{r,s,p+1}. 
\end{eqnarray*}
Hence we find that, up to a constant multiple, 
$Q_{2,2v}$ is given by
\begin{eqnarray*}
 Q_{2,2v} &=& \sum_{r+s+p=v} \frac{1}{r!s!p!}
                             (d_1/2-3/2+v)_{v-r}\  
                             (d_2/2-3/2+v)_{v-s} \\
 &&\qquad \qquad  (-((d_1+d_2)/2-3/2+v))_{v-p}\    
                  (P^*_0)^r (P_2^*)^s (P^*_1)^p
\end{eqnarray*}
where we have, again, used $(x)_n$ for $\prod_{i=0}^{n-1} (x-i)$.

Since it is obvious that 
\begin{eqnarray*}
&& \D^p( \D^r F(Z)\ \D^s G(Z) ) = \\
  && \left(
   (P^*_1(\partial_{Z_1},\partial_{Z_2}))^p\
   (P^*_0(\partial_{Z_1},\partial_{Z_2}))^r\
   (P^*_2(\partial_{Z_1},\partial_{Z_2}))^s  F(Z_1) G(Z_2)
   \right) \vert_{Z=Z_1=Z_2}
\end{eqnarray*}
where $\D = \det(\partial_{Z})$ we obtain exactly the formula in 
Theorem 1.2  of {\it loc. cit.}.

\section{A vector valued generalisation}
\label{vector}

In this section we describe vector valued differential 
operators $D$ such that  
$$
D(F(Z_{1})G(Z_{2}))_{|Z_{1}=Z_{2}}
$$
is a vector valued modular form if $F$ and $G$ are modular forms.

Let us introduce some notation first.
For any representation $\rho$ of $\GL{n}$ we denote by $d(\rho)$ its 
dimension. Furthermore, for any even positive integers $d_{1}$ and 
$d_{2}$,  we denote by $\Hq_{n,\rho}(d_{1},d_{2})$ 
the space of $d(\rho)$ dimensional vectors of 
$O(d_{1}) \times O(d_{2})$-invariant 
pluri-harmonic polynomials 
$P(X,X^{'})=(P_{i}(X,X^{'}))_{1 \leq i \leq d(\rho)}$ 
such that $P(A X,A X^{'}) = \rho(A)P(X,X^{'})$.

The main result of ref.\ \cite{Ibukiyama} not only includes 
the case of invariant pluri-harmonic polynomials in
$\Hq_{n,v}(d_{1},d_{2})$ ({\it c.f.} Theorem \ref{Ibuthm})
but also applies to the case of  pluri-harmonic polynomials 
in $\Hq_{n,\rho}(d_{1},d_{2})$. It is shown in {\it loc. cit.}
that the vector valued differential operators corresponding 
to the latter polynomials map two modular forms to a vector 
valued modular form. We will, therefore, concentrate on giving 
some examples of pluri-harmonic polynomials in 
$\Hq_{n,\rho}(d_{1},d_{2})$.

More precisely, we will consider only the representations 
of $\GL{n}$ which are of the form $\rho_{m,v}=\det^{v} \Sym^{m}$
where $\Sym^{m}$ is the symmetric tensor representation of degree $m$. 
When $n = 2$ these representations exhaust all polynomial 
representations of $\GL{n}$. 

For given $m$ and $v$ the Young diagram corresponding to 
$\rho_{m,v}$ is given by  
$(v+m,v,...,v)$ where the depth is $n$ if $v \neq 0$ and $1$ if $v = 0$.
For the sake of simplicity, we assume from now on that $d_{1} \geq 2n$ and 
$d_{2} \geq 2n$.
Then one can easily see that $\dim(\Hq_{n,\rho_{m,v}}(d_{1},d_{2})) = 1$ 
if and only if $m$ and $v$ are even.   
We want to give explicit bases of (some of) the  one dimensional 
spaces $\Hq_{n,\rho_{m,v}}(d_1,d_2)$. 

Firstly, we consider the simplest case, {\it i.e.} the case of the 
symmetric representation $\rho_{m,0}$. 
In this case everything can be reduced to the genus one case 
which was already discussed in \S\ref{genus1}.

Let $u_{1},...,u_{n}$ be $n$ independent variables and, 
for any multi-index $\nu=(\nu_{1},...,\nu_{n})\in\N_0^n$, 
write $u^{\nu}$ for $\prod_{i=1}^{n}u_{i}^{\nu_{i}}$ 
and  $|\nu|$ for $\nu_{1} + \cdots + \nu_{n}$. 
Denote by $I(m)$ the set of all multi-indices with $|\nu|=m$
and set $S = \sum_{i,j=1}^{n}R_{ij}u_{i}u_{j}$  
and $S^{'} = \sum_{i,j=1}^{n}R^{'}_{ij}u_{i}u_{j}$.  
\begin{prop}
\label{vectorprop1}
Let $m$ be even and let $Q(r,r^{'})$ be a basis of $\Hq_{1,m}(d_1,d_2)$.
For each multi-index $\nu$ with $|\nu| = m$ we define $Q_{\nu}(R,R^{'})$ 
by $Q(S,S^{'}) = \sum_{|\nu|=m}Q_{\nu}(R,R^{'})u^{\nu}$. 
Then the vector $(Q_{\nu}(R,R^{'}))_{\nu \in I(m)}$ 
gives a basis of the one dimensional space $\Hq_{n,\rho_{m,0}}(d_{1},d_{2})$.
\end{prop}
{\it Proof.} 
Since it is clear that 
\[
(Q_{\nu}(ARA^{t},AR^{'}A^{t}))_{\nu \in I(m)} 
=
\rho_{m,0}(A)(Q_{\nu}(R,R^{'}))_{\nu \in I(m)}
\]
all we should do is to prove that the polynomials $Q_{\nu}(R,R^{'})$
are harmonic.
This is easily proved as follows. Using that $Q(r,r^{'})$ is harmonic 
one obtains 
\[
\left(2d_{1}\frac{\partial Q(r,r^{'})}{\partial r} 
+ 4r\frac{\partial^{2} Q(r,r^{'})}{\partial r^{2}}\right)
+\left(2d_{2}\frac{\partial Q(r,r^{'})}{\partial r^{'}} 
+ 4r^{'}\frac{\partial^{2} Q(r,r^{'})}{\partial r^{'2}}\right)
= 0.
\]
Furthermore, a simple calculation shows that 
\begin{eqnarray*}
L^{(d_{1})}_{ij}(Q(S,S^{'}))
& = & 
u_{i}u_{j}
(2d_{1}\frac{\partial Q}{\partial r} 
+ 4S\frac{\partial^{2} Q}{\partial r^{2}})\, ,
\\
L^{(d_{2})}_{ij}(Q(S,S^{'}))
& = & 
u_{i}u_{j}
(2d_{2}\frac{\partial Q}{\partial r^{'}} 
+ 4S^{'}\frac{\partial^{2} Q}{\partial r^{'2}})
\end{eqnarray*}
so that the proposition becomes obvious.
\qed

Finally, we give some examples for the case of representations
of mixed type when $n=2$. For $(m+v,v) = (m+2,2)$ with even $m$ 
an invariant pluri-harmonic polynomial in $\Hq_{2,\rho_{m,2}}$
is given by $(Q_{\nu}(R,R^{'}))_{\nu \in I(m)}$ where 
$Q_{\nu}(R,R^{'})$ is the coefficients of $u^{\nu}$ of the following 
polynomial $Q(S,S^{'})$.
\begin{eqnarray*}
Q(S,S^{'}) & = & 
Q_{2,d_{1},d_{2}}(R_{1},R_{2})F_{m,d_{1}+2,d_{2}+2}(S,S^{'})
\\ 
&&+\frac{1}{2}((d_{2}-1)P_{0}S-(d_{1}-1)P_{2}S^{'})
(\frac{\partial F_{m.d_{1}+2,d_{2}+2}}{\partial r}
-\frac{\partial F_{m.d_{1}+2,d_{2}+2}}{\partial s})(S,S^{'}). 
\end{eqnarray*}
Here we have denoted by $Q_{v,d_{1},d_{2}}$ the non-zero invariant
pluri-harmonic polynomial in $\Hq_{2,v}(d_{1},d_{2})$ normalised 
as in section \ref{genus2} and by $F_{m,d_{1},d_{2}}(r,s)$ 
a non-zero polynomial in $\Hq_{1,m}(d_{1},d_{2})$. 

When $(m+2,2) = (4,2)$, for example, then $Q$ this is given by 
\begin{eqnarray*}
Q(S,S^{'}) 
& = & \ \ 
(d_{2}-1)d_{2}(d_{2}+2)P_{0}S - (d_{1}-1)(d_{2}-1)(d_{2}+2)P_{1}S
\\
 & &
+(d_{1}-1)(d_{1}+2)(d_{4}+4)P_{2}S
- (d_{1}+4)(d_{2}-1)(d_{2}+2)P_{0}S^{'}
\\ & &
+(d_{1}-1)(d_{1}+2)(d_{2}-1)P_{1}S^{'} 
- d_{1}(d_{1}-1)(d_{1}+2)P_{2}S^{'}.
\end{eqnarray*}

\section{Conclusion}
\label{conclusion}

In this paper we have described certain spaces of invariant pluri-harmonic 
polynomials. These polynomials are in one-to-one correspondence with
Rankin-Cohen type differential operators which, in the case of 
non-constant polynomials, map two Siegel modular forms to a cusp form.
In particular, we have derived a set of recursion 
equations (Theorem \ref{vgeneral}) which uniquely (up to multiplication 
by non-zero elements in $\C^*$) determine an invariant 
pluri-harmonic polynomial $Q_{n,v}$ in $\Hq_{n,v}(d_1,d_2)$ (for
$d_1\ge n$ and $d_2\ge n$). 
Although the recursion equations can easily be solved for
any {\it numerical} values of $n$ and $v$ we have not been able 
to give the closed explicit formulas for the solutions for general 
$n$ and $v$. However, in several examples we have obtained such 
formulas for the solutions.
In addition, we have discussed certain vector valued bilinear 
differential operators.

Let us conclude with a few remarks and point out some interesting 
open questions in connection with our results.

Firstly, the polynomials $Q_{n,v}\in \Hq_{n,v}(d_1,d_2)$ are in 
one-to-one correspondence with the differential operators $D$ 
defined in Theorem \ref{Ibuthm} only if $d_1\ge n$ and $d_2\ge n$.
Throughout our analysis we have always assumed that 
this condition is satisfied.
For fixed $n$ and $v$ the polynomial $Q_{n,v}$ depends only on 
$d_1$ and $d_2$ and, with a suitable normalisation, is well 
defined and non-vanishing even if $d_1< n$ or $d_2< n$.
Therefore, one might speculate that the differential operators 
corresponding to $Q_{n,v}$ in the latter cases also map any two automorphic 
forms to an automorphic form (for $n=1,2$ this follows from the results 
in ref.\ \cite{Cohen,ChoieEholzer}). 
In contrast to the situation for $d_1\ge n$ and $d_2\ge n$ 
the dimension of $\Hq_{n,v}(d_1,d_2)$ can be larger than one
if $d_1$ or $d_2$ are `small'.  For example if $d_1=d_2=n-1$ 
and $v=2$ then $\Hq_{n,2}(n-1,n-1)$ is spanned by $P_0$ and
$P_n$ (note that, in this case, the pluri-harmonic polynomials 
$P_0(X X^t, X' {X'}^t)$ and $P_n(X X^t, X' {X'}^t)$ are identically zero).
In this case the differential operators corresponding to $P_0$ and
$P_n$ satisfy the commutation relation ({\it c.f.} \S\ref{basics}).
For odd $n$, however, all modular forms of weight $(n-1)/2$ on 
congruence subgroups are singular so that these differential operators
act as zero (this special case is already contained in the results of 
Kapitel III, \S6 of ref.\ \cite{Freitag} and ref.\ \cite{Resnikoff}).
Furthermore, note that if $d_1$ and $d_2$ are such that 
$d_1+d_2<2n$ then there is no pluri-harmonic polynomial 
in $\P_{n,v}(d_1,d_2)$ ({\it c.f.} the discussion in ref. 
\cite{KashiwaraVergne}).
The last remarks show that it would be interesting to compute 
the dimension of $\Hq_{n,v}(d_1,d_2)$ for $d_1+d_2 \ge 2n$
and either $d_1$ or $d_2$ small and to understand the relation between 
covariant differential operators and invariant pluri-harmonic 
polynomials in this case.

Secondly, the  differential operators $D$ give rise to 
differential operators for Jacobi forms of higher degree. 
The relation between the Rankin-Cohen type operators for $n=2$ 
and the corresponding differential operators for Jacobi forms has 
been discussed in detail in ref.\ \cite{ChoieEholzer}. 
One would expect that --like in the case $n=2$-- the dimension of the
space of covariant differential operators for higher degree Jacobi 
forms is, for fixed $v$, generically greater than one. 
It would be interesting if one could compute this dimension and 
obtain, as in the case of $n=2$, closed explicit formulas.

Thirdly, in ref.\ \cite{Zagier} certain algebraic structures --called 
Rankin-Cohen algebras-- have been defined using only the Rankin-Cohen 
operators for $n=1$. There is an obvious generalisation of this
definition for the case of arbitrary $n$ using the corresponding 
Rankin-Cohen type differential operators studied in this paper. 
Therefore, one is naturally led to the question whether one can
describe the structure of these generalised Rankin-Cohen algebras in 
an independent way analogous to the case $n=1$ ({\it cf.} the Theorem 
in {\it loc. cit.}). Furthermore, it would be interesting to know if 
other examples of these algebraic structures can be found in 
mathematical nature.

Fourthly, it seems quite natural to look for applications of the 
Rankin-Cohen type operators constructed in this paper which involve
theta series and/or Eisenstein series.

Finally, let us mention that there are other very interesting types of 
pluri-harmonic polynomials. The polynomials considered in this paper 
correspond to the second case in ref.\ \cite{Ibukiyama}.
The first case considered in {\it loc. cit.} is concerned with differential 
operators $D$ acting on automorphic functions $F$ on $\H_n$ such that 
the restriction of $D(F)$ to $\H_{n_1}\times\cdots\times\H_{n_r}$
with $n=n_1+\cdots+n_r$ is again an automorphic form 
and the corresponding pluri-harmonic polynomials \cite{IbukiyamaZagier}.

\medskip\noindent
{\bf Acknowledgements}

We would like to thank D.\ Zagier for many inspiring discussions and, 
in particular, for the explicit solution of the recursion equation
in \S\ref{v4}.

W.E. would also like to thank Y.\ Choie for many stimulating discussions.

The second author is grateful to the Max-Planck-Institut f\"ur
Mathematik in Bonn for its hospitality during his visit when this 
work was done.


\end{document}